\documentclass[11pt,a4paper]{article}
\usepackage{caption2}%---change the caption of figures
\usepackage{CJK}
\usepackage{tabls}
\usepackage{bm}
\usepackage{multirow}
\usepackage{amssymb}
\usepackage{amsfonts}
\usepackage{mathrsfs}
\usepackage{empheq}
\usepackage{amsmath}
\usepackage{amsthm}
\usepackage{graphicx}
\usepackage{color}
\usepackage{indentfirst}
\usepackage{hyperref}%pour PDFLaTex
\usepackage{float}
\usepackage{cases}
\usepackage[titletoc]{appendix}
\usepackage{bm}

\usepackage[colorinlistoftodos,english]{todonotes}
\allowdisplaybreaks%

\def\e{\eqref}

\def\i1n{i=1,\cdots,n}
\def\j1n{j=1,\cdots,n}
\def\ij1n{i,j=1,\cdots,n}

\def \i{\mathrm i}
\def\e{\varepsilon}
\def\R{\mathbf R}

 \numberwithin{equation}{section}
\theoremstyle{definition}

 \newtheorem{thm}{\indent Theorem}[section]
 
 \newtheorem{lem}{\indent Lemma}[section]

 \newtheorem{rem}{\indent Remark}[section]

\theoremstyle{definition}
\theoremstyle{theorem}
\theoremstyle{lemma}

\newcommand{\be}{\begin{equation}}
\newcommand{\ee}{\end{equation}}
\newcommand{\beq}{\begin{equation*}}
\newcommand{\eeq}{\end{equation*}}

\begin{document}

\begin{CJK*}{GB}{gbsn}
\title{Blow-up and lifespan estimate to a nonlinear wave equation in Schwarzschild spacetime}

\author{Ning-An Lai \\ School of Mathematics Sciences, Zhejiang Normal University, Jinhua 321004, P. R. China \thanks{ninganlai@zjnu.edu.cn. Corresponding author}
\and Yi Zhou \\ School of Mathematics Sciences, Fudan University, Shanghai 200433, P. R. China \thanks{yizhou@fudan.edu.cn }}

\pagestyle{myheadings} \markboth{Blow-up for semilinear wave equation in Schwarzschild spacetime}{Blow-up for semilinear wave equation in Schwarzschild spacetime}
\maketitle

\begin{abstract}

We study the semilinear wave equation with power type nonlinearity and small initial data in Schwarzschild spacetime. If the nonlinear exponent $p$ satisfies $2\le p<1+\sqrt 2$, we establish the sharp upper bound of lifespan estimate, while for the most delicate critical power $p=1+\sqrt2$, we show that the lifespan satisfies\\
 \[
 T(\e)\le \exp\left(C\e^{-(2+\sqrt 2)}\right),
 \]
 the optimality of which remains to be proved. The key novelty is that the compact support of the initial data can be close to the event horizon. By combining the global existence result for $p>1+\sqrt 2$ obtained by Lindblad et al.(Math. Ann. 2014), we then give a positive answer to the interesting question posed by Dafermos and Rodnianski(J. Math. Pures Appl. 2005, the end of the first paragraph in page $1151$): $p=1+\sqrt 2$ is exactly the critical power of $p$ separating stability and
blow-up.
\end{abstract}

\textbf{Keywords}: Semilinear wave equations; blow-up; Schwarzschild spacetime; lifespan; event horizon

\textbf{2010 MSC}: 35L70, 58J45

\section{Introduction}

Nonlinear wave equations in Schwarzschild spacetime attract more and more attention, since it is natural to generalize the classical results for nonlinear wave equations in flat Minkowski spactime to the black hole spacetime. Schwarzschild metric is the first analytic solution to the vacuum Einstein equation, which was derived by Schwarzschild in 1915. And according to Birkhoff¡¯s theorem it is also a
unique spherical symmetry solution of the vacuum Einstein equation. The explicit expression of the metric is
\begin{equation}\label{Smetric}
\begin{aligned}
g_S=F(r)dt^2-F(r)^{-1}dr^2-r^2d\omega^2,
\end{aligned}
\end{equation}
where $F(r)=1-\frac{2M}{r}$, and $M$ is the Newtonian mass. Noting that asymptotically $r\rightarrow\infty$ or $M\rightarrow 0$ the Schwarzschild metric reduces to the
Minkowski metric
\[
ds^2=dt^2-dr^2-r^2d\omega^2.
\]

The study of classical Cauchy problem to semilinear wave equations in Minkowski spacetime
\begin{equation}
\label{Mmain}
\left\{
\begin{aligned}
& u_{tt} - \Delta u = |u|^p, \quad\text{in $[0, T)\times\R^n$}, \\
& u(x,0)=\e f(x), \quad u_t(x,0)=\e g(x), \quad x\in\R^n
\end{aligned}
\right.
\end{equation}
has a long time history, and it has now been determined there exists a critical power $p_c(n)>1, n\ge 2$, which solves the quadratic equation
\[
(n-1)p^2-(n+1)p-2=0.
\]
Such kind problem is known as Strauss conjecture: for $1<p\le p_c(n)$, the solution will blow up in a finite time, while for $p>p_c(n)$ the solution exists globally in time, see \cite{Georgive, Glassey1, Glassey2, John, LZ, Schaeffer, Sideris, Yordanov, Zhou2} and references therein. It is easy to see that $p_c(3)=1+\sqrt 2$. If there is no global solution, it is then also interesting to estimate the lifespan($T(\e)$) with respect to the small parameter $\e$. We are now clear that
there exist two positive constants $c$ and $C$
such that the lifespan satisfies for $n\ge 2$ and
$\max(1,2/(n-1))<p<p_c(n)$
\begin{equation}\label{1c}
\begin{aligned}
c\varepsilon^{\frac{-2p(p-1)}{\gamma(n, p)}}\leq T(\varepsilon)\leq C\varepsilon^{\frac{-2p(p-1)}{\gamma(n, p)}},
\end{aligned}
\end{equation}
where $\gamma(n, p)=2+(n+1)p-(n-1)p^2>0$. For $(n, p)=(2, 2)$,
\begin{equation}
\left\{
\begin{array}{ll}
\exists \lim\limits_{\e\rightarrow 0^{+}} a(\e)^{-1}T(\e)>0, & \mathrm{if}~\int_{\R^2}g(x)dx\neq 0,\\
\exists \lim\limits_{\e\rightarrow 0^{+}} \e T(\e)>0,& \mathrm{if}~\int_{\R^2}g(x)dx= 0,
\end{array}
\right.
\end{equation}
where $a(\e)$ denotes a number satisfying
\[
a^2\e^2\log(1+a)=1.
\]
For $1<p<2$ and $n=2$,
\begin{equation}
%\label{np22}
\left\{
\begin{array}{ll}
c\e^{-\frac{p-1}{3-p}}\le T(\e)\le C\e^{-\frac{p-1}{3-p}},&\mathrm{if}~\int_{\R^2}g(x)dx\neq 0,\\
c\e^{-\frac{2p(p-1)}{\gamma(2, p)}}\le T(\e)\le C\e^{-\frac{2p(p-1)}{\gamma(2, p)}},&\mathrm{if} ~\int_{\R^2}g(x)dx= 0.
\end{array}
\right.
\end{equation}
For the critical case $(p=p_c(n), n\ge 2)$,
\begin{equation}\label{1d}
\begin{aligned}
\exp(c\varepsilon^{-p(p-1)})\leq T(\varepsilon)\leq \exp(C\varepsilon^{-p(p-1)}).
\end{aligned}
\end{equation}
See \cite{Takamura6, LZ, LZbook, Lindblad3, Lindblad2, Takamura2, Takamura, Zhou3, Zhou5} and and the introduction in \cite{LLWW}.

It is natural to consider the corresponding Cauchy problem in Schwarzschild spacetime
\begin{equation}
\label{Smain}
\left\{
\begin{aligned}
& \Box_{g_S} u = |u|^p, \quad (t, x)\in \R_t^+\times \Sigma, \\
& u(x,0)=\e f(x), \quad u_t(x,0)=\e g(x), \quad x\in \Sigma,
\end{aligned}
\right.
\end{equation}
where $g_S$ denotes the Schwarzschild metric presented in \eqref{Smetric} and $\R_t^+\times \Sigma$ is called the
exterior of the black hole:
\begin{equation}\nonumber\\
\begin{aligned}
\R_t^+\times \Sigma=\R_t^+\times (2M, \infty)\times \mathbb{S}^2.\\
\end{aligned}
\end{equation}
The D'Alembert operator associated with the Schwarzschild metric $g$ can be written explicitly as
\begin{equation}\label{3}
\begin{aligned}
\Box_{g_S}=\frac{1}{F(r)}\Big(\partial_t^2-\frac{F}{r^2}\partial_r(r^2F)\partial_r-\frac{F}{r^2}\Delta_{\mathbb{S}^2}\Big),
\end{aligned}
\end{equation}
where $\Delta_{\mathbb{S}^2}$ denotes the standard Laplace-Beltrami operator on $\mathbb{S}^2$. At first glance, a divergent coefficient(when $r$ approaches to $2M$) appears in the operator, which may cause essential difficulty near the event horizon($r=2M$). The study of power type nonlinear wave equations in Schwarzschild(-like) and Kerr black hole spacetime was initiated by Nicolas \cite{Nico1, Nico2}, in which the global existence results were established for
\[
\Box_gu+m^2u+\lambda |u|^2u=0,~~~\lambda>0
\]
with large initial data, and $g$ denotes the Schwarzschild-like or Kerr metric. After that, more and more attention is paid to the small data Cauchy problem. Catania and Georgiev \cite{Catania} first studied the blow-up phenomenon for \eqref{Smain}. They considered radial solution in the Regge-Wheeler coordinate
\begin{equation}\label{RH}
r^*=r+2M\ln(r-2M),
\end{equation}
and proved that the solution will blow up for $1<p<1+\sqrt2$, by choosing special initial data
\[
f(r^*)=\e \chi_0(r^*-r_0^*(\e)),~~g(r^*)=\e \chi_1(r^*-r_0^*(\e)),
\]
where $\chi_j\in C_0^{\infty}(\R)(j=1, 2)$ satisfy
\[
\begin{aligned}
&\chi_j(r^*)\ge 0, r^*\in \R,\\
&\chi_j(r^*)=1,  r^*\in [-R/2, R/2],\\
&supp~\chi_j\subset [-R, R]
\end{aligned}
\]
for some positive constant. And
\[
r_0^*(\e)=\e^{-\theta}
\]
for some positive constant $\theta$ depending on $p$. Noting that $r^*$ varies from $-\infty$ to $+\infty$ as $r$ varies from $2M$ to $+\infty$, it means that the chosen data above have compact support far away from the event horizon, and the distance depends on the small parameter, see the picture below.
\begin{center}

 \hspace{-4mm}
~~~~~~\includegraphics[scale=0.5]{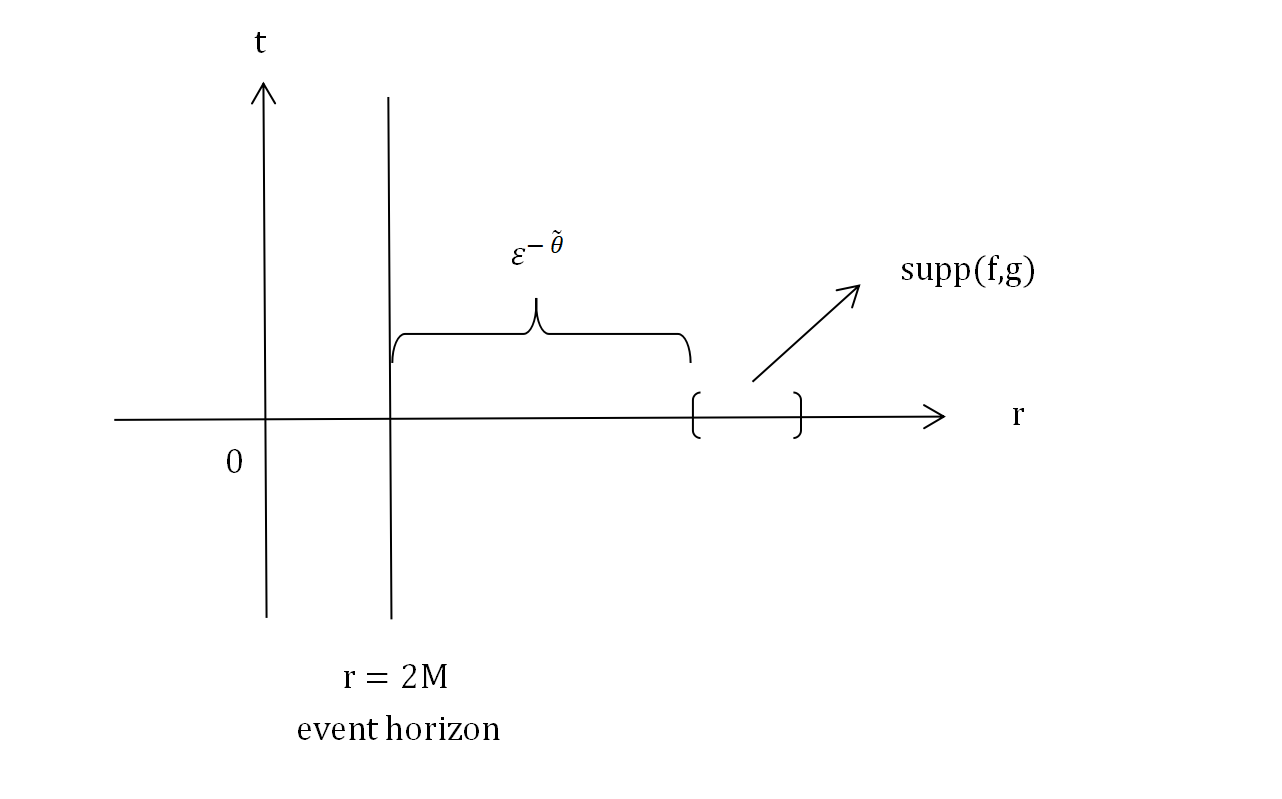}\\

\end{center}

Dafermos and Rodnianski \cite{DR} first studied the global existence for the corresponding small data problem, in which they showed the radial solution for $p>4$ in Schwarzschild or more generally a Reissner-Nordstr\"{o}m spacetime exists globally in time. Blue and Sterbenz \cite{BS} generalized the global existence of general solution for \eqref{Smain} to $p>3$. Marzuola et al \cite{MMT} proved the global-in-time Strichartz estimates for wave equations in Schwarzschild spacetime and furthermore established global existence for \eqref{Smain} with $p=5$. Lindblad et al \cite{LMS} then established global existence result for $p>1+\sqrt2$, which also holds for Kerr black hole with small angular momentum if the small data have compact support. This result takes a giant step forward for the question posed by Dafermos and Rodnianski \cite{DR}, and it was generalized by Metcalfe-Wang \cite{MW} to the slowly rotating Kerr spacetime with non-compact supported small data. Wang \cite{Wan1} established the following lower bound of the lifespan estimate
\begin{equation}
\label{lowerbound}
T(\varepsilon)\ge\left\{
\begin{aligned}
& C\varepsilon^{-\frac{p(p-1)}{1+2p- p^2}}, \quad &2\le p<1+\sqrt2, \\
& \exp\left(C\varepsilon^{-2\sqrt2}\right),  & p=1+\sqrt2,\\
\end{aligned}
\right.
\end{equation}
which also holds for the corresponding problem posed on Kerr spacetime with small angular momentum($a\ll M$), see \cite{DFW}. The lower bound for $2\le p<1+\sqrt 2$ above coincides with that in \eqref{1c} by setting $n=3$, while the existence time for the \lq\lq critical" power $p=1+\sqrt 2$ seems to have a long gap, comparing to the one in \eqref{1d}. Very recently, Lin et al. \cite{LLM} established blow-up result for $p\in [\frac{3}{2} , 2]$, without assuming that the supports of the initial data should be far away from the black hole. We also refer the reader to \cite{Wan} for more introduction to the Cauchy problem \eqref{Smain}.

One should also notice that in \cite{Luk} Luk first studied the semilinear wave equations with derivative nonlinearity in slowly
rotating Kerr spacetimes$(a\ll M)$, and showed global existence for small data solution assuming the nonlinear term satisfies the null condition. Without null condition assumption, the authors \cite{LaiZ} proved that the semilinear problem in Schwarzschild spacetime has no global solution.

The main target of this paper is to show that $p=1+\sqrt 2$ is indeed the critical power to \eqref{Smain}, thus, to show there is no global solution for $2\le p\le1+\sqrt 2$, without any assumption on the distance between the event horizon and the support of the initial data. What is more, we will prove the sharpness of the lifespan estimate for $2\le p<1+\sqrt 2$, while for the critical power $p=1+\sqrt2$, a possibly sharp lifespan estimate from above will be established, since there is a gap comparing to the second lower bound in \eqref{lowerbound}. If we try to get blow-up result for small initial data problem, it always means that the nonlinear term will dominate the linear effect, this is why we may succeed in showing blow-up for relative small power($1<p\le p_c(n)$) in the flat Minkowski spacetime. However, for our concerned problem \eqref{Smain}, it is easy to see from the wave operator \eqref{3} with Schwarzschild metric that a factor asymptotically approaching to $0$ as coming close to the event horizon($r=2M$) will appear. This factor will make competition with the nonlinear term and it may prevent finite time blow-up in some sense. This fact is also the reason why the assumption that the support of the initial data should be far away from the event horizon is needed for small data problem in \cite{Catania}. We in this paper divide the radial variable into two parts: $r^*<\widetilde{C}_0$ and $r^*\ge \widetilde{C}_0$ for some positive constant $\widetilde{C}_0$, which corresponds to $2M<r<2M+C_0$ and $r\ge 2M+C_0$ for some other positive constant $C_0$ respectively. For the latter part, the wave will not touch the event horizon then all the estimates are almost the same as that of the corresponding flat Minkowski case, hence the nonlinear term will have polynomial growth at most. We call this area polynomial zone, which is labelled as $Z_{pol}$ in the picture below. For the former part $2M<r<2M+C_0$($r^*<\widetilde{C}_0$), the upper bound of the nonlinear term will increase exponentially with respect to $t$, which seems impossible to obtain a lifespan estimate of polynomial type. This exponential zone is labelled by $Z_{exp}$ in the picture below. However, if we use the test function
\[
\psi_{\lambda}(t, r)=e^{-\lambda t}\phi_{\lambda},
\]
where $\phi_{\lambda}$ solves for $\lambda>0$
\begin{equation}\label{test}
\begin{aligned}
\frac{1}{r^2}\partial_r\big(r(r-2M)\partial_r\phi_{\lambda}\big)=
\frac{\lambda^2}{1-\frac{2M}{r}}\phi_{\lambda},
\end{aligned}
\end{equation}
we then can gain some exponential decay in $Z_{exp}$ by choosing cut-off functions delicately, and hence the exponential increase mentioned above can be cancelled out. Inspired by this observation, we find that the effect for blow-up caused by $Z_{pol}$ area will dominate $Z_{exp}$ zone. We then use a family of solutions to the elliptic equation \eqref{test} to construct a powerful test function $b_q$(see Lemma \ref{bq}) for the critical power $p=1+\sqrt2$, and furthermore to establish blow-up result and lifespan estimate for the solution involved in $Z_{pol}$ zone.
\begin{center}

 \hspace{-4mm}
~~~~~~\includegraphics[scale=0.5]{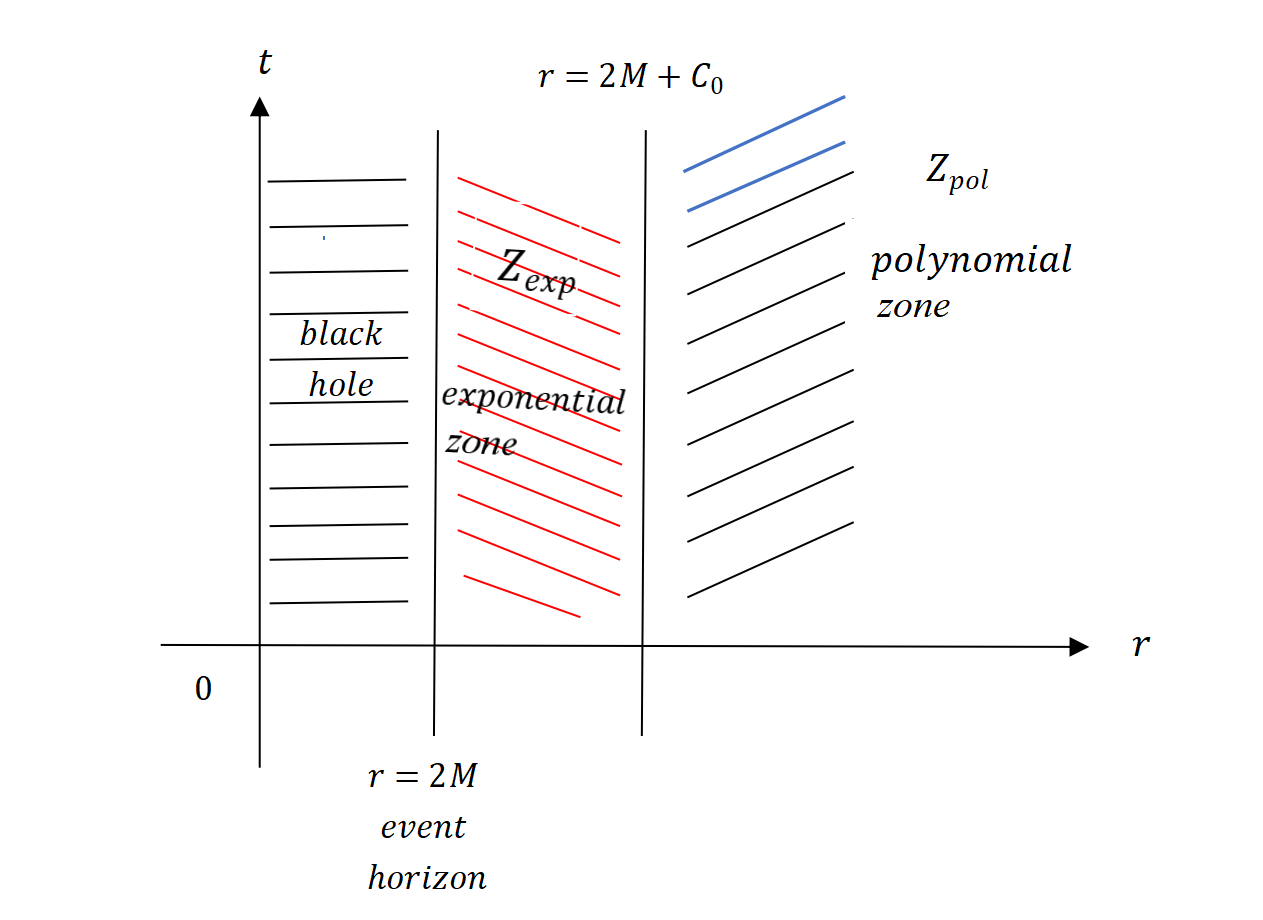}\\

\end{center}

We state our main results as
\begin{thm}\label{thm1}
Let $2\leq p<1+\sqrt 2$. Suppose the initial data in \eqref{Smain} are nonnegative and do not vanish identically, and the supports of the data satisfy
\begin{equation}\label{supp}
\begin{aligned}
supp~f(x), g(x)\subset \{2M+R_1\le r\le 2M+R_2\}\times\mathbb{S}^2
\end{aligned}
\end{equation}
for some positive constants $R_1<R_2$. Then the solution of problem \eqref{Smain} blows up in a finite time. Furthermore, the upper bound of lifespan estimate satisfies
 \begin{equation}\label{lf}
\begin{aligned}
T(\varepsilon)\leq C\varepsilon^{-\frac{p(p-1)}{1+2p-p^2}}.
\end{aligned}
\end{equation}
Hereafter, $C$ denotes a generic positive constant independent of $\e$, and the values may change from line to line.
\end{thm}

\begin{rem}
The assumption \eqref{supp} posed on the initial data in our results can be depicted as
\begin{center}

 \hspace{-4mm}
~~~~~~\includegraphics[scale=0.5]{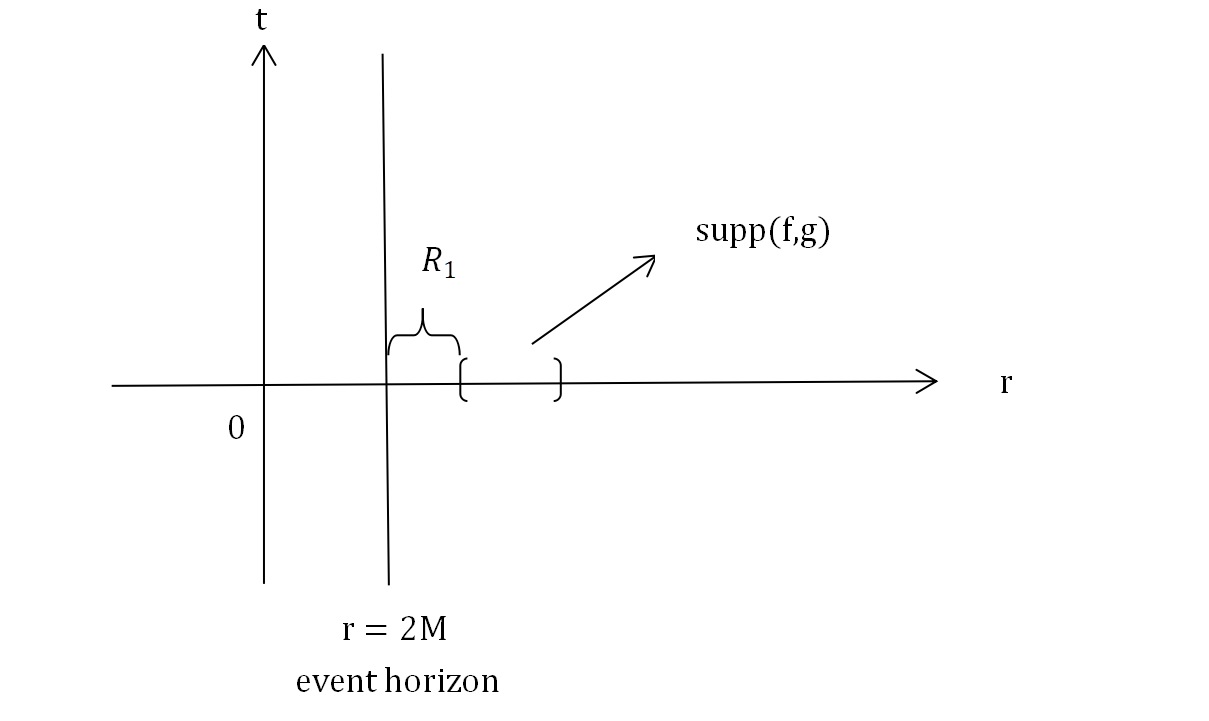}\\

\end{center}
\end{rem}

\begin{rem}
In $r^*$ coordinate, the support assumption \eqref{supp} implies that there exists some positive constant $R$ such that
\begin{equation}\label{supp1}
\begin{aligned}
supp~f(x), g(x)\subset \{|r^*|\le R\}\times\mathbb{S}^2.
\end{aligned}
\end{equation}
\end{rem}

\begin{thm}\label{thm2}
Let $p=1+\sqrt 2$. Suppose the initial data satisfy the same conditions as in Theorem \ref{thm1}. Then the solution of problem \eqref{Smain} blows up in a finite time. Furthermore, the upper bound of the lifespan estimate satisfies
 \begin{equation}\label{lf2}
\begin{aligned}
T(\varepsilon)\leq \exp\left(C\e^{-p(p-1)}\right)=\exp\left(C\e^{-(2+\sqrt 2)}\right).
\end{aligned}
\end{equation}
\end{thm}

\begin{rem}

In some sense, $p=1+\sqrt 2$ is the \lq\lq only" power left to be solved, according to the question posed at the end of the first paragraph in page $1151$ in \cite{DR}. The above theorem shows that it also belongs to the blow-up case. Actually, we also first show the finite time blow-up result for $2\le p<1+\sqrt2$ when the supports of the initial data are close to the black hole.

\end{rem}

\begin{rem}

Compared to the lower bound of lifespan for $p=1+\sqrt2$ in \eqref{lowerbound}, there still exists a gap between the existence and the blow-up time.
\begin{center}

 \hspace{-4mm}
~~~~~~\includegraphics[scale=0.25]{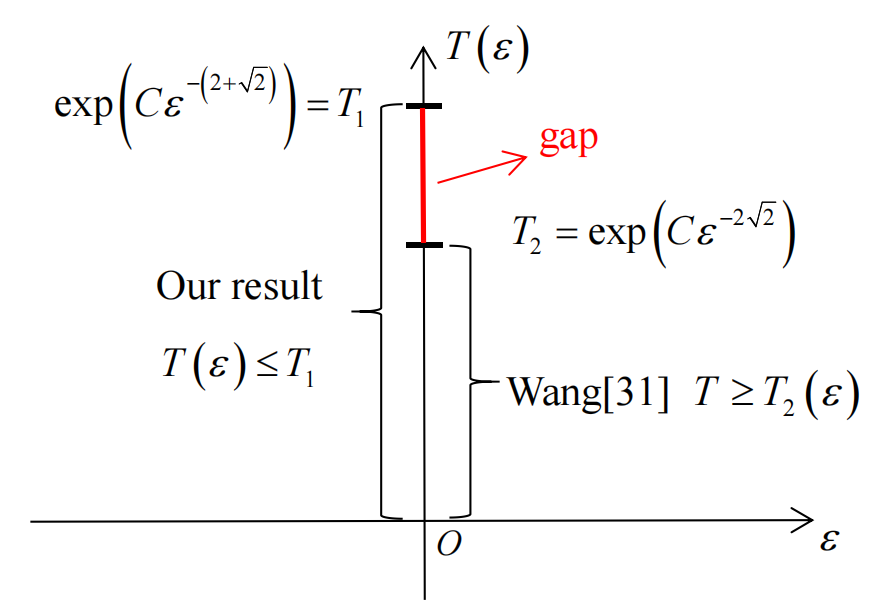}\\

\end{center}

\end{rem}

We organize the paper as follows. In Section 2, we list the ingredients used to construct the appropriate test function. In Section 3, the upper bound of the spacetime integral of the nonlinear term far away from the event horizon is established, by using the space and time cut-off functions, while in Section 4 the upper bound of the spacetime integral of the nonlinear term in the whole exterior of the black hole is obtained, by using only a time cut-off function. In Section 5, we demonstrate the proof for Theorem 1.1, which corresponds to the upper bound of lifespan estimate for subcritical powers. The main idea is to obtain the optimal lower bounds for the space time integrals of the nonlinear term, then the result follows by combining the upper bounds established in Section 3 and 4. In Section 6, we show the proof for Theorem 1.4, which corresponds to the blow-up and upper bound of lifespan estimate for the critical power $p=1+\sqrt2$. The key ingredient is that we use a family of test functions with parameter $\lambda$ to construct a new test function with better decay.

\section{Preliminaries}

For blow-up result, the normal way is to find some appropriate test functions, and then to construct some functional including the solution which approaches to $\infty$ at some time. And hence the upper bound of lifespan estimate follows. Besides, we choose some smooth cut-off functions delicately
\[
\eta(t)=
 \left\{
 \begin{array}{ll}
 1 &  0\le t\le\frac13,\\
% decreasing & \mbox{for}\ \frac12<t<1,\\
 0 &  t\ge \frac23,\\
 \end{array}
 \right.
\eta_T(t)=\eta\left(\frac tT\right),
\]
\[
\alpha(r^*)=
 \left\{
 \begin{array}{ll}
 0 &  -\infty\le r^*\le\frac18,\\
% decreasing & \mbox{for}\ \frac12<t<1,\\
 1 &  r^*\ge \frac14,\\
 \end{array}
 \right.
\alpha_T(r^*)=\alpha\left(\frac {r^*}T\right)
\]
and
\[
\chi(\theta)=
 \left\{
 \begin{array}{ll}
 1 &  0\le \theta\le\frac34,\\
% decreasing & \mbox{for}\ \frac12<t<1,\\
 0 &  \theta\ge \frac56,\\
 \end{array}
 \right.
\chi_T(\theta)=\chi\left(\frac {\theta}T\right),
\]
where
\begin{equation}\label{T}
T\in [T_0, T(\e)),~~~T_0=\max\{8(4M+e), 12(R+R_3), 16R\},
\end{equation}
with $R, R_3$ the positive constants appearing in \eqref{supp1} and \eqref{Philam} respectively. The next lemma concerns to our key test function.
\begin{lem}\label{lphi1}
The elliptic equation with $0<\lambda\le\lambda_0$ for some fixed $\lambda_0(>\max\{\frac{8}{Mp(p-1)}, 1\})$
\begin{equation}\label{test1}
\begin{aligned}
\frac{1}{r^2}\partial_r\big(r(r-2M)\partial_r\phi_{\lambda}\big)=
\frac{\lambda^2}{1-\frac{2M}{r}}\phi_{\lambda}
\end{aligned}
\end{equation}
admits a family of nonnegative solutions $\phi_{\lambda}(r)$ satisfy
\begin{equation}\label{phi1}
\begin{aligned}
\phi_{\lambda}(r) \thicksim \frac 1re^{\lambda r^*}\thicksim\frac 1{(1+\lambda)r}e^{\lambda r^*}
\end{aligned}
\end{equation}
and
\begin{equation}\label{dphi1}
\begin{aligned}
|\partial_r\phi_{\lambda}(r)|\lesssim\frac {\lambda}{r-2M}e^{\lambda r^*}.
\end{aligned}
\end{equation}
\end{lem}
\textit{\emph Proof.} In Lemma 5.4 in \cite{Catania}, it has been proved that the equation
\[
\partial_{r^*r^*}\varphi_{\lambda}-\frac{2M\left(1-\frac{2M}{r}\right)}{r^3}
\varphi_{\lambda}=\lambda^2\varphi_{\lambda}
\]
admits a family of nonnegative solutions $\varphi_{\lambda}(r^*)$ satisfying
\begin{equation}\label{var1}
\varphi_{\lambda}(r^*)\thicksim e^{\lambda r^*}.
\end{equation}
Set
\[
\phi_{\lambda}(r)=\frac{\varphi_{\lambda}}{r},
\]
then direct computation implies that $\phi_{\lambda}(r)$ solves \eqref{test1}, and the asymptotic behavior \eqref{phi1} follows from \eqref{var1}.

For \eqref{dphi1}, equation \eqref{test1} yields
\begin{equation}\label{estdphi}
\begin{aligned}
\partial_r\phi_{\lambda}&=\frac{\lambda^2}{r(r-2M)}\int_{2M}^{r}\frac{\tau^3}
{\tau-2M}\phi_{\lambda}d\tau\\
&=\frac{\lambda^2}{r(r-2M)}\int_{-\infty}^{r^*}\tau^2\phi_{\lambda}ds.\\
\end{aligned}
\end{equation}
We next divide the estimate into two parts. For $r^*\le 0$, we have $\tau, r\thicksim 1$, then plugging \eqref{phi1} into \eqref{estdphi} yields
\begin{equation}\label{estdphi1}
\begin{aligned}
|\partial_r\phi_{\lambda}|&\lesssim\frac{\lambda^2}{r-2M}
\int_{-\infty}^{r^*}e^{\lambda s}ds\\
&\lesssim \frac{\lambda}{r-2M}e^{\lambda r^*}.
\end{aligned}
\end{equation}
For $r^*>0$, we have $r\thicksim 1+r^*$ and then
\begin{equation}\label{estdphi2}
\begin{aligned}
|\partial_r\phi_{\lambda}|&\lesssim\frac{\lambda^2}{r(r-2M)}\left(
\int_{-\infty}^{0}e^{\lambda s}ds+\int_0^{r^*}se^{\lambda s}ds\right)\\
&\lesssim\frac{\lambda^2}{r(r-2M)}\left(\frac{1}{\lambda}+\frac{1}{\lambda}r^*
e^{\lambda r^*}\right)\\
&\lesssim \frac{\lambda}{r-2M}e^{\lambda r^*},
\end{aligned}
\end{equation}
hence \eqref{dphi1} follows from \eqref{estdphi1} and \eqref{estdphi2}.
\section{Upper Bound For The Nonlinear Term Far Away From Event Horizon}

Multiplying the equation in \eqref{Smain} with $\eta_T^{2p'}(t)
\alpha_T^{2p'}(r^*)r^2$, and then integrating with respect to $\omega, r, t$, we get
\begin{equation}\label{s31}
\begin{aligned}
&\int_0^T\int_{2M}^{\infty}\int_{\mathbb{S}^2}\left[\frac{r^2}{1-\frac{2M}{r}}
\partial_t^2u-\partial_r\big(r(r-2M)\partial_ru\big)-\Delta_{\mathbb{S}^2}u\right]
\eta_{T}^{2p'}\alpha_T^{2p'}d\omega drdt\\
=&\int_0^T\int_{2M}^{\infty}\int_{\mathbb{S}^2}|u|^p\eta_T^{2p'}
\alpha_T^{2p'}r^2d\omega drdt\\
\triangleq &F_0(T).\\
\end{aligned}
\end{equation}
Due to the assumption \eqref{supp1}, the support of the solution satisfies
\begin{equation}\label{supp3}
supp u~\subset \Big\{r^*\Big| |r^*|\le t+R\Big\}\times \mathbb{S}^2.
\end{equation}
Hence by integration by part we get from \eqref{s31}
\begin{equation}\label{s32}
\begin{aligned}
F_0(T)\lesssim&\int_0^T\int_{2M}^{\infty}
\int_{\mathbb{S}^2}\Big|\frac{r^2}{1-\frac{2M}{r}}u\partial_t^2\eta_T^{2p'}\alpha
_T^{2p'}\Big|d\omega drdt\\
&+\int_0^T\int_{2M}^{\infty}
\int_{\mathbb{S}^2}\big|u\eta_T^{2p'}(2r-2M)\partial_r\alpha_T^{2p'}\big|d\omega drdt\\
&+\int_0^T\int_{2M}^{\infty}
\int_{\mathbb{S}^2}\big|u\eta_T^{2p'}(r^2-2Mr)\partial_r^2\alpha_T^{2p'}\big|d\omega drdt\\
\triangleq& I_1+I_2+I_3.\\
\end{aligned}
\end{equation}
Since
\[
|\partial_t^2\eta_T^{2p'}|\lesssim T^{-2}\eta_{T}^{2p'-2},
\]
and if $r^*\ge \frac T8\ge 4M+e$, then
\[
r\ge 2M+e,~~r^*\thicksim r,
\]
and hence
\[
\frac{r}{r-2M}\thicksim C,
\]
we may estimate $I_1$ by H\"{o}lder inequality
\begin{equation}\label{I1}
\begin{aligned}
&\int_0^T\int_{2M}^{\infty}
\int_{\mathbb{S}^2}\Big|\frac{r^2}{1-\frac{2M}{r}}u\partial_t^2\eta_T^{2p'}\alpha
_T^{2p'}\Big|d\omega drdt\\
\lesssim &T^{-2}F_0^{\frac1p}(T)\left(\int_{\frac T3}^{\frac{2T}{3}}
\int_{r^*\ge \frac T8}\int_{\mathbb{S}^2}r^2\big(\frac{r}{r-2M}\big)^{p'}d\omega drdt\right)^{\frac 1{p'}}\\
\lesssim &T^{-2}F_0^{\frac1p}(T)\left(\int_{\frac T3}^{\frac{2T}{3}}
\int_{4M+e}^{t+R}\int_{\mathbb{S}^2}r^2\big(\frac{r}{r-2M}\big)^{p'}\frac
{r-2M}{r}d\omega dr^*dt\right)^{\frac 1{p'}}\\
\lesssim &T^{-2}F_0^{\frac1p}(T)\left(\int_{\frac T3}^{\frac{2T}{3}}
\int_{4M+e}^{t+R}(r^*)^2 dr^*dt\right)^{\frac 1{p'}}\\
\lesssim &T^{\frac 4{p'}-2}F_0^{\frac1p}(T).
\end{aligned}
\end{equation}
For $r^*\ge 4M+e$, it has
\[
r-M\thicksim r\thicksim r^*,
\]
then $I_2$ can be estimated as
\begin{equation}\label{I2}
\begin{aligned}
&\int_0^T\int_{2M}^{\infty}
\int_{\mathbb{S}^2}\big|u\eta_T^{2p'}(2r-2M)\partial_r\alpha_T^{2p'}\big|d\omega drdt\\
\lesssim &T^{-1}F_0^{\frac1p}(T)\left(\int_{0}^{T}
\int_{r^*\ge \frac T8}\int_{\mathbb{S}^2}r^{-\frac{2}{p-1}}r^{p'}d\omega drdt\right)^{\frac 1{p'}}\\
\lesssim &T^{-1}F_0^{\frac1p}(T)\left(\int_{0}^{T}
\int_{4M+e}^{t+R}(r^{*})^{\frac{p-2}{p-1}}d\omega dr^*dt\right)^{\frac 1{p'}}\\
\lesssim &T^{\frac 4{p'}-2}F_0^{\frac1p}(T),
\end{aligned}
\end{equation}
and we can get the similar estimate for $I_3$
\begin{equation}\label{I3}
\begin{aligned}
\int_0^T\int_{2M}^{\infty}
\int_{\mathbb{S}^2}\big|u\eta_T^{2p'}(r^2-2Mr)\partial_r^2\alpha_T^{2p'}\big|d\omega drdt
\lesssim T^{\frac 4{p'}-2}F_0^{\frac1p}(T).
\end{aligned}
\end{equation}
Since it also holds
\[
r-2M\thicksim r\thicksim r^*,~~for~r^*>4M+e,
\]
we then get the upper bound for $F_0(T)$ by combining
\eqref{s32}, \eqref{I1}, \eqref{I2} and \eqref{I3}
\begin{equation}\label{F0}
\begin{aligned}
F_0(T)
\lesssim T^{4-2p'}.
\end{aligned}
\end{equation}

\section{Upper Bound For The Nonlinear Term Outside The Black Hole}

Comparing to the estimate of $F_0(T)$, we have no space cut-off for $F_1(T)$, so we have to complete the estimate by dividing the radial space variable into two cases: $r^*\ge 4M+e$ and $r^*\le 4M+e$. This time we multiply the equation in \eqref{Smain} with $\eta_T^{2p'}r^2$ to get
\begin{equation}\label{s41}
\begin{aligned}
&\int_0^T\int_{2M}^{\infty}\int_{\mathbb{S}^2}\left[\frac{r^2}{1-\frac{2M}{r}}
\partial_t^2u-\partial_r\big(r(r-2M)\partial_ru\big)-\Delta_{\mathbb{S}^2}u\right]
\eta_{T}^{2p'}d\omega drdt\\
=&\int_0^T\int_{2M}^{\infty}\int_{\mathbb{S}^2}|u|^p\eta_T^{2p'}
r^2d\omega drdt\\
\triangleq &F_1(T).\\
\end{aligned}
\end{equation}
Again by integration by parts and H\"{o}lder inequality, we have
\begin{equation}\label{s42}
\begin{aligned}
F_1(T)\le& \int_0^T\int_{2M}^{\infty}\int_{\mathbb{S}^2}\frac{r^2}{1-\frac{2M}{r}}u
\partial_t^2\eta_T^{2p'}d\omega drdt\\
\lesssim&T^{-2}F_1^{\frac1p}(T)\left(\int_{\frac T3}^{\frac{2T}{3}}
\int_{|r^*|\le t+R}\int_{\mathbb{S}^2}r^2\left(\frac{r}{r-2M}\right)^{p'}d\omega drdt\right)^{\frac{1}{p'}}\\
\lesssim& T^{-2}F_1^{\frac1p}(T)\Bigg[\int_{\frac T3}^{\frac{2T}{3}}
\int_{4M+e}^{ t+R}r^2\left(\frac{r}{r-2M}\right)^{p'}\frac{r-2M}{r} dr^*dt\\
&+\int_{\frac T3}^{\frac{2T}{3}}
\int_{-t-R}^{4M+e}r^2\left(\frac{r}{r-2M}\right)^{p'}\frac{r-2M}{r} dr^*dt\Bigg]^{\frac1{p'}}\\
\triangleq&T^{-2}F_1^{\frac1p}(T)(I_4+I_5)^{\frac1{p'}}.
\end{aligned}
\end{equation}
We may control $I_4$ in the similar way as that of $I_1$, thus
\begin{equation}\label{I4}
\begin{aligned}
I_4\lesssim T^{4}.\\
\end{aligned}
\end{equation}
But for $I_5$, since
\[
r^*\le 4M+e\Leftrightarrow 2M\le r\le 2M+e,
\]
and hence
\begin{equation}\label{ln}
\begin{aligned}
r-2M\sim e^{\frac{r^*-C}{2M}}\thicksim e^{\frac{r^*}{2M}},
\end{aligned}
\end{equation}
we then have
\begin{equation}\label{I5}
\begin{aligned}
I_5=&\int_{\frac T3}^{\frac{2T}{3}}
\int_{-t-R}^{4M+e}r^2\left(\frac{r}{r-2M}\right)^{p'}\frac{r-2M}{r} dr^*dt\\
\lesssim&\int_{\frac T3}^{\frac{2T}{3}}
\int_{-t-R}^{4M+e}e^{-\frac{1}{p-1}\frac{r^*}{2M}}dr^*dt\\
\lesssim&\int_{\frac T3}^{\frac{2T}{3}}e^{\frac{t}{2M(p-1)}}dt\\
\lesssim&e^{\frac{T}{3M(p-1)}}.\\
\end{aligned}
\end{equation}
By combining \eqref{s42}, \eqref{I4} and \eqref{I5} we come to
\begin{equation}\label{F1}
\begin{aligned}
F_1\lesssim& T^{4-2p'}+T^{-2p'}e^{\frac{T}{3M(p-1)}}\\
\lesssim& T^{-2p'}e^{\frac{T}{3M(p-1)}}.\\
\end{aligned}
\end{equation}

\section{Proof For Theorem \ref{thm1}}

In order to prove Theorem \ref{thm1}, we use
\begin{equation}\label{Philam}
\Phi_{\lambda}(t, r)=\eta_T^{2p'}(t)\chi_T^{2p'}(t-r^*+R_3)e^{-\lambda t}\phi_{\lambda}(r)r^2
\end{equation}
as the test function, where $R_3$ is a positive constant independent of $\e$. Also, we can choose $R_3$ and $T$(large enough) such that
\begin{equation}\label{supp2}
\begin{aligned}
&supp~\{\chi_T(R_3-r^*)=1\}\cap supp~(f, g)=supp~(f, g),\\
&supp~\{\partial_t\chi_T(R_3-r^*)\}\cap supp~(f, g)=\varnothing.
\end{aligned}
\end{equation}
Multiplying the equation in \eqref{Smain} with $\Phi_{\lambda}(t, r)$ and making integration by parts, one has
\begin{equation}\label{s51}
\begin{aligned}
&\int_0^T\int_{2M}^{\infty}\int_{\mathbb{S}^2}\partial_t
\left(\frac{r^2}{1-\frac{2M}{r}}
\partial_tu\eta_T^{2p'}\chi_T^{2p'}e^{-\lambda t}\phi_{\lambda}\right)d\omega drdt\\
&-\int_0^T\int_{2M}^{\infty}\int_{\mathbb{S}^2}\partial_t
\left[\frac{r^2}{1-\frac{2M}{r}}
u(\eta_T^{2p'}\chi_T^{2p'}e^{-\lambda t}\phi_{\lambda})_t\right]d\omega drdt\\
&+\int_0^T\int_{2M}^{\infty}\int_{\mathbb{S}^2}
\frac{r^2}{1-\frac{2M}{r}}
u\partial_t^2(\eta_T^{2p'}\chi_T^{2p'}e^{-\lambda t}\phi_{\lambda})d\omega drdt\\
&-\int_0^T\int_{2M}^{\infty}\int_{\mathbb{S}^2}
u\partial_r\left[r(r-2M)\partial_r(\eta_T^{2p'}\chi_T^{2p'}e^{-\lambda t}\phi_{\lambda})\right]d\omega drdt\\
=&\int_0^T\int_{2M}^{\infty}\int_{\mathbb{S}^2}|u|^p\eta_T^{2p'}
\chi_T^{2p'}e^{-\lambda t}\phi_{\lambda}r^2d\omega drdt,\\
\end{aligned}
\end{equation}
which implies according to the support assumption \eqref{supp2} that
\begin{equation}\label{s52}
\begin{aligned}
&C_1\e\\
\le&\int_0^T\int_{2M}^{\infty}\int_{\mathbb{S}^2}
\frac{r^2}{1-\frac{2M}{r}}
u\Big(\partial_t^2(\eta_T^{2p'}\chi_T^{2p'})e^{-\lambda t}\phi_{\lambda}-2\lambda\partial_t\eta_T^{2p'}\chi_T^{2p'}
e^{-\lambda t}\phi_{\lambda}\Big)d\omega drdt\\
&-2\lambda\int_0^T\int_{2M}^{\infty}\int_{\mathbb{S}^2}
\frac{r^2}{1-\frac{2M}{r}}
u\eta_T^{2p'}\partial_t\chi_T^{2p'}
e^{-\lambda t}\phi_{\lambda}d\omega drdt\\
&-\int_0^T\int_{2M}^{\infty}\int_{\mathbb{S}^2}
u\Big(2r(r-2M)\eta_T^{2p'}\partial_r\chi_T^{2p'}e^{-\lambda t}\partial_r\phi_{\lambda}\\
&+r(r-2M)\eta_T^{2p'}\partial_{rr}\chi_T^{2p'}e^{-\lambda t}\phi_{\lambda}+2(r-M)\eta_T^{2p'}\partial_r\chi_T^{2p'}e^{-\lambda t}\phi_{\lambda}\Big)d\omega drdt\\
\triangleq&I_6+I_7+I_8+I_9+I_{10}+I_{11},\\
\end{aligned}
\end{equation}
where
\[
C_1=\int_{2M}^{\infty}\int_{\mathbb{S}^2}\frac{r^2}{1-\frac{2M}{r}}
\big(\lambda f(x)+g(x)\big)
\phi_{\lambda}(r)d\omega drdt.
\]
For the terms $I_6-I_{11}$, the key point is that each of them contains at least one derivative over one of the cut-off functions, thanks to the equation \eqref{test1} in Lemma \ref{lphi1}. We estimate $I_7$ first.
\begin{equation}\label{I7}
\begin{aligned}
I_7\triangleq
&-2\lambda\int_0^T\int_{2M}^{\infty}\int_{\mathbb{S}^2}
\frac{r^2}{1-\frac{2M}{r}}
u\partial_t\eta_T^{2p'}\chi_T^{2p'}
e^{-\lambda t}\phi_{\lambda}d\omega drdt\\
\lesssim& T^{-1}\Big[\int_0^T\int_{2M}^{\infty}\int_{\mathbb{S}^2}
\frac{r^2}{1-\frac{2M}{r}}\big|u\eta_{T}^{2p'-1}\eta_T{'}\chi_T^{2p'}\alpha_T^{2p'}
\big|e^{-\lambda t}\phi_{\lambda}d\omega drdt\\
&+ \int_0^T\int_{2M}^{\infty}\int_{\mathbb{S}^2}
\frac{r^2}{1-\frac{2M}{r}}\big|u\eta_{T}^{2p'-1}\eta_T{'}\chi_T^{2p'}
(1-\alpha_T^{2p'})
\big|e^{-\lambda t}\phi_{\lambda}d\omega drdt\Big]\\
\triangleq &T^{-1}({I_7}_1+{I_7}_2),\\
\end{aligned}
\end{equation}
where $I_{71}$ can be estimated by $F_0(T)$
\begin{equation}\label{I71}
\begin{aligned}
I_{71}\triangleq&
\int_0^T\int_{2M}^{\infty}\int_{\mathbb{S}^2}
\frac{r^2}{1-\frac{2M}{r}}\big|u\eta_{T}^{2p'-1}\eta_T{'}\chi_T^{2p'}\alpha_T^{2p'}
\big|e^{-\lambda t}\phi_{\lambda}d\omega drdt\\
\lesssim& F_0(T)^{\frac1p}\left(\int_{\frac T3}^{\frac{2T}{3}}\int_{4M+e}^{t+R}
r^2\left(1-\frac{2M}{r}\right)^{1-p'}r^{-p'}e^{-\lambda p'(t-r^*)}dr^*dt\right)^{\frac1{p'}}\\
\lesssim& F_0(T)^{\frac1p}\left(\int_{\frac T3}^{\frac{2T}{3}}\int_{4M+e}^{t+R}
(r^*)^{2-p'}e^{-\lambda p'(t-r^*)}dr^*dt\right)^{\frac1{p'}}\\
\lesssim& F_0(T)^{\frac1p}(T^{3-p'})^{\frac1{p'}},
\end{aligned}
\end{equation}
where an elementary integral estimate \eqref{16a} below has been used.
\begin{lem}[\cite{LaiZ}, Lemma 3.3]\label{lem6}
Given any $\alpha\geq 0, \beta>0$ and $L>0$, there exists a positive constant $C$ such that
\begin{equation}\label{16a}
\begin{aligned}
\int_{0<s\leq t+L}(1+s)^{\alpha}e^{-\beta(t-s)}ds\leq C(t+L)^{\alpha}.
\end{aligned}
\end{equation}
\end{lem}
\textit{\emph Proof.} For convenience, we outline the proof. Dividing the integral in \eqref{16a} into two parts
\begin{equation}\label{16b}
\begin{aligned}
&\int_{0<s\leq t+L}(1+s)^{\alpha}e^{-\beta(t-s)}ds\\
=&\int_{0}^{\frac{t+L}{2}}(1+s)^{\alpha}e^{-\beta(t-s)}ds+\int_{\frac{t+L}{2}}^{t+L}
(1+s)^{\alpha}e^{-\beta(t-s)}ds\\
\leq &Ce^{-\beta\cdot \frac{t-L}{2}}\int_{0}^{\frac{t+L}{2}}(1+s)^{\alpha}ds+C(t+L)^{\alpha}\int_{\frac{t+L}{2}}^{t+L}
e^{-\beta(t-s)}ds\\
\leq &C(t+L)^{\alpha},\\
\end{aligned}
\end{equation}
which is the desired inequality \eqref{16a}.

We control $I_{72}$ by dividing the integral into two parts with the assumption $\lambda\ge \frac{1}{2Mp}$
\begin{equation}\label{I72}
\begin{aligned}
I_{72}\triangleq&
\int_{\frac T3}^T\int_{2M}^{\infty}\int_{\mathbb{S}^2}
\frac{r^2}{1-\frac{2M}{r}}\big|u\eta_{T}^{2p'-1}\eta_T{'}\chi_T^{2p'}
(1-\alpha_T^{2p'})
\big|e^{-\lambda t}\phi_{\lambda}d\omega drdt\\
\lesssim& F_1(T)^{\frac1p}\Bigg[\int_{\frac T3}^{\frac{2T}{3}}\int_{4M+e}^{\frac T4}
r^2\left(1-\frac{2M}{r}\right)^{1-p'}r^{-p'}e^{-\lambda p'(t-r^*)}dr^*dt\\
&+\int_{\frac T3}^{\frac{2T}{3}}\int_{-t-1}^{4M+e}
r^2\left(1-\frac{2M}{r}\right)^{1-p'}r^{-p'}e^{-\lambda p'(t-r^*)}dr^*dt\Bigg]^{\frac1{p'}}\\
\lesssim& F_1(T)^{\frac1p}\Bigg[\int_{\frac T3}^{\frac{2T}{3}}\int_{4M+e}^{\frac T4}
(r^*)^{2-p'}r^{-p'}e^{-\lambda p'(t-r^*)}dr^*dt\\
&+\int_{\frac T3}^{\frac{2T}{3}}\int_{-t-R}^{4M+e}
e^{\frac{r^*}{2M}(1-p')}e^{-\lambda p'(t-r^*)}dr^*dt\Bigg]^{\frac1{p'}}\\
\lesssim& F_1(T)^{\frac1p}\Bigg[\int_{\frac T3}^{\frac{2T}{3}}T^{2-p'}e^{-\lambda p't}dt\int_{4M+e}^{\frac T4}
e^{\lambda p'r^*}dr^*\\
&+\int_{\frac T3}^{\frac{2T}{3}}e^{-\lambda p't}dt\int_{-t-R}^{4M+1}
e^{\frac{r^*}{2M(p-1)}(2M\lambda p-1)}dr^*\Bigg]^{\frac1{p'}}\\
\lesssim& F_1(T)^{\frac1p}\left(T^{3-p'}e^{-\frac{\lambda p'T}{12}}+T^2e^{-\frac{\lambda p'T}{3}}\right)^{\frac1{p'}}.
\end{aligned}
\end{equation}
We next estimate $I_9$. Since in the support of $\partial_r\chi_T(t-r^*+R_3)$ one has
\[
\frac{3T}{4}\le t-r^*+R_3\le \frac{5T}{6},
\]
which implies $t\ge \frac T3$ by combining \eqref{T} and \eqref{supp3}, and hence
\begin{equation}\label{I9}
\begin{aligned}
I_9\triangleq
&-\int_0^T\int_{2M}^{\infty}\int_{\mathbb{S}^2}
u2r(r-2M)\eta_T^{2p'}\partial_r\chi_T^{2p'}e^{-\lambda t}\partial_r\phi_{\lambda}d\omega drdt\\
\lesssim& T^{-1}\Big[\int_0^T\int_{2M}^{\infty}\int_{\mathbb{S}^2}
\big|ur(r-2M)\eta_{T}^{2p'}\chi_T^{2p'-1}\partial_r\chi_T\alpha_T^{2p'}
\big|e^{-\lambda t}\partial_r\phi_{\lambda}d\omega drdt\\
&+ \int_0^T\int_{2M}^{\infty}\int_{\mathbb{S}^2}
\big|ur(r-2M)\eta_{T}^{2p'}\chi_T^{2p'-1}\partial_r\chi_T(1-
\alpha_T^{2p'})
\big|e^{-\lambda t}\partial_r\phi_{\lambda}d\omega drdt\Big]\\
\triangleq &T^{-1}({I_9}_1+{I_9}_2),\\
\end{aligned}
\end{equation}
where $I_{91}$ can be controlled by $F_0(T)$
\begin{equation}\label{I91}
\begin{aligned}
I_{91}\triangleq
&\int_{\frac T3}^T\int_{2M}^{\infty}\int_{\mathbb{S}^2}
\big|ur(r-2M)\eta_{T}^{2p'}\chi_T^{2p'-1}\partial_r\chi_T\alpha_T^{2p'}
\big|e^{-\lambda t}\partial_r\phi_{\lambda}d\omega drdt\\
\lesssim&\int_{\frac T3}^T\int_{2M}^{\infty}\int_{\mathbb{S}^2}
\big|ur^2\eta_{T}^{2p'}\chi_T^{2p'-1}\alpha_T^{2p'}
\big|e^{-\lambda t}\partial_r\phi_{\lambda}d\omega drdt\\
\lesssim& F_0^{\frac1p}(T)\left(\int_{\frac T3}^{\frac{2T}{3}}\int_{4M+e}^{t+R}r^2(r-2M)^{-p'}e^{-\lambda p'(t-r^*)}\frac{r-2M}{r}dr^*dt\right)^{\frac{1}{p'}}\\
\lesssim& F_0^{\frac1p}(T)\left(\int_{\frac T3}^{\frac{2T}{3}}\int_{4M+e}^{t+R}(r^*)^{2-p'}e^{-\lambda p'(t-r^*)}dr^*dt\right)^{\frac{1}{p'}}\\
\lesssim& F_0^{\frac1p}(T)(T^{3-p'})^{\frac{1}{p'}},\\
\end{aligned}
\end{equation}
while $I_{92}$ can be estimated by $F_1(T)$ as \eqref{I72}
\begin{equation}\label{I92}
\begin{aligned}
I_{92}\triangleq
&\int_0^T\int_{2M}^{\infty}\int_{\mathbb{S}^2}
\big|ur(r-2M)\eta_{T}^{2p'}\chi_T^{2p'-1}\partial_r\chi_T(1-
\alpha_T^{2p'})
\big|e^{-\lambda t}\partial_r\phi_{\lambda}d\omega drdt\\
\lesssim&\int_0^T\int_{2M}^{\infty}\int_{\mathbb{S}^2}
\big|ur^2\eta_{T}^{2p'}\chi_T^{2p'-1}\alpha_T^{2p'}
\big|e^{-\lambda t}\partial_r\phi_{\lambda}d\omega drdt\\
\lesssim& F_1^{\frac1p}(T)\Big(\int_{\frac T3}^{\frac{2T}{3}}\int_{4M+e}^{\frac T4}r^2(r-2M)^{-p'}e^{-\lambda p'(t-r^*)}\frac{r-2M}{r}dr^*dt\\
&+\int_{\frac T3}^{\frac{2T}{3}}\int_{-t-R}^{4M+e}r^2(r-2M)^{-p'}e^{-\lambda p'(t-r^*)}\frac{r-2M}{r}dr^*dt\Big)^{\frac{1}{p'}}\\
\lesssim& F_1(T)^{\frac1p}\Big(\int_{\frac T3}^{\frac{2T}{3}}T^{2-p'}e^{-\lambda p't}dt\int_{4M+e}^{\frac T4}
e^{\lambda p'r^*}dr^*\\
&+\int_{\frac T3}^{\frac{2T}{3}}e^{-\lambda p't}dt\int_{-t-R}^{4M+1}
e^{\frac{r^*}{2M(p-1)}(2M\lambda p-1)}dr^*\Big)^{\frac1{p'}}\\
\lesssim& F_1(T)^{\frac1p}\left(T^{3-p'}e^{-\frac{\lambda p'T}{12}}+T^2e^{-\frac{\lambda p'T}{3}}\right)^{\frac1{p'}}.
\end{aligned}
\end{equation}
It is easy to see that $I_6$ and $I_8$ can be estimated in the similar way as that of $I_7$, while the terms $I_{10}$ and $I_{11}$ can be done in the way as that of $I_9$, and hence we finally come to by combining \eqref{F0} and \eqref{F1}
\begin{equation}\label{final1}
\begin{aligned}
\e\lesssim&F_0^{\frac1p}T^{-1}(T^{3-p'})^{\frac{1}{p'}}+F_1^{\frac1p}T^{-1}
\left(T^{3-p'}e^{-\frac{\lambda p'T}{12}}+T^2e^{-\frac{\lambda p'T}{3}}\right)^{\frac1{p'}}\\
\lesssim&T^{\frac{p^2-2p-1}{p(p-1)}}+T^{-\frac{2p'}{p}+\frac{3}{p'}-1}
e^{T\left(\frac{1}{3Mp(p-1)}-\frac{\lambda}{12}\right)}\\
\lesssim&T^{\frac{p^2-2p-1}{p(p-1)}},\\
\end{aligned}
\end{equation}
if we choose $\lambda$ large enough such that
\begin{equation}\label{final}
\lambda>\frac{4}{Mp(p-1)}.
\end{equation}
The lifespan estimate \eqref{lf} comes from \eqref{final1}.

\section{Proof For Theorem \ref{thm2}}

Set
\[
\lambda=\frac{8}{Mp(p-1)}
\]
in the proof in last section, then from the first inequality in \eqref{final1}, there exists a positive constant $C_1$ depending on $p, M$ such that
\begin{equation}\label{60}
\begin{aligned}
C_1\e&\le F_0^{\frac1p}T^{-1}(T^{3-p'})^{\frac{1}{p'}}+F_1^{\frac1p}T^{-2
+\frac{3}{p'}}e^{-\frac{\lambda T}{12}}\\
&\le F_0^{\frac1p}T^{-1}(T^{3-p'})^{\frac{1}{p'}}+F_1^{\frac1p}T^{-2
+\frac{3}{p'}}e^{-\frac{T}{3Mp(p-1)}}.\\
\end{aligned}
\end{equation}
We claim that the second term in the right hand side of \eqref{60} should be less than $\frac{C_1\e}{2}$. If not, we than have by combining \eqref{F1}
\begin{equation}\label{61}
\begin{aligned}
\frac{C_1\e}{2}&\lesssim T^{4-4\sqrt 2},
\end{aligned}
\end{equation}
where we plug into the value $p=1+\sqrt 2$. This fact implies that the lifespan for $p=1+\sqrt 2$ is at least of polynomial type, which contradict the exponential lower bound established by \cite{Wan1}, see \eqref{lowerbound} above. This claim further implies that
\begin{equation}\label{62}
\begin{aligned}
F_0^{\frac1p}T^{-1}(T^{3-p'})^{\frac{1}{p'}}\ge \frac{C_1\e}{2},
\end{aligned}
\end{equation}
which yields
\begin{equation}\label{63}
\begin{aligned}
F_0\gtrsim \e^pT^{3-p}.
\end{aligned}
\end{equation}
Actually, from the process to get \eqref{final1}, we may get a more precise lower bound than \eqref{63}, thus
\begin{equation}\label{64}
\begin{aligned}
\overline{F}_0\triangleq\int_{\frac T3}^{\frac{2T}{3}}\int_{2M}^{\infty}\int_{\mathbb{S}^2}|u|^p\eta_T^{2p'}
\alpha_T^{2p'}r^2d\omega drdt\gtrsim \e^pT^{3-p},
\end{aligned}
\end{equation}
which definitely implies
\begin{equation}\label{65}
\begin{aligned}
\widetilde{F}_0\triangleq\int_{\frac T{16}}^{\frac{2T}{3}}\int_{2M}^{\infty}\int_{\mathbb{S}^2}|u|^p\eta_T^{2p'}
\alpha_T^{2p'}r^2d\omega drdt\gtrsim \e^pT^{3-p}.
\end{aligned}
\end{equation}

In the following, we are going to use a similar method employed in \cite{LT1, LLWW}.
\begin{lem}\label{bq}
Let $\phi_{\lambda}$ be the one in Lemma \ref{lphi1} and
\[
q=1-\frac1p=1-\frac{1}{1+\sqrt2}=2-\sqrt 2.
\]
Set
\[
b_q(t, r)=\int_0^1e^{-\lambda t}\phi_{\lambda}(r)\lambda^{q-1}d\lambda,
\]
then it holds
\begin{equation}\label{066}
\begin{aligned}
\partial_tb_q=b_{q+1},~~\frac{r^2}{1-\frac{2M}{r}}\partial_t^2b_q-\partial_r
\left(r(r-2M)\partial_rb_q\right)=0,\\
\end{aligned}
\end{equation}
and for $r^*\ge4M+e$
\begin{equation}\label{66}
\begin{aligned}
&|\partial_rb_q|\lesssim b_{q+1},\\
&b_q\thicksim (t+R)^{-q},\\
&b_{q+1}\lesssim (t+R)^{-1}(t+R+1-r^*)^{-q}.\\
\end{aligned}
\end{equation}
\end{lem}
\textit{\emph Proof.} The two identities in $\eqref{066}$ can be obtained by direct computations, noting that $\phi_{\lambda}$ satisfies \eqref{test1}. For $\eqref{66}_1$, by \eqref{dphi1} and \eqref{phi1}, one has for $r^*\ge 4M+e$
\begin{equation}\label{67}
\begin{aligned}
|\partial_rb_q|&\le \int_0^1e^{-\lambda t}|\partial_r\phi_{\lambda}|\lambda^{q-1}d\lambda\\
&\lesssim \int_0^1e^{-\lambda t}\frac{\lambda}{r-2M}e^{\lambda r^*}\lambda^{q-1}d\lambda\\
&\lesssim  \int_0^1e^{-\lambda t}\frac{r}{r-2M}\phi_{\lambda}\lambda^{q}d\lambda\\
&\lesssim b_{q+1}.\\
\end{aligned}
\end{equation}
We then show
\begin{equation}\label{68}
b_q\gtrsim (t+R)^{-q},
\end{equation}
for which the key ingredient is the uniform positive lower bound for $\phi_\lambda$. Noting \eqref{estdphi} and the nonnegative of $\phi_{\lambda}$, we know $\partial_r\phi_\lambda$ is nonnegative and hence $\phi_{\lambda}$ is nondecreasing with respect to $r$ for each fixed $\lambda$, and hence for $r^*\ge 4M+e(r\ge 2M+e)$ we have for $\lambda\in(0, 1]$
\begin{equation}\label{69}
\begin{aligned}
\phi_\lambda(r)&\ge \phi_{\lambda}(2M+e)\\
&\gtrsim \frac{e^{\lambda(4M+e)}}{2M+e}\\
&\gtrsim \frac1{2M+e}.
\end{aligned}
\end{equation}
With this in hand, we then get
\begin{equation}\label{601}
\begin{aligned}
b_q(t, r)&\gtrsim \int_{\frac{1}{2(t+R)}}^{\frac{1}{t+R}}e^{-\lambda t}\phi_\lambda \lambda^{q-1}d\lambda\\
&\gtrsim \int_{\frac{1}{2(t+R)}}^{\frac{1}{t+R}}e^{-\lambda(t+R)} \lambda^{q-1}d\lambda\\
&=(t+R)^{-q}\int_{\frac12}^{1}e^{-\theta}\theta^{q-1}d\theta\\
&\gtrsim (t+R)^{-q}.\\
\end{aligned}
\end{equation}
For the upper bound of $b_q, b_{q+1}$, we divide the proof into two cases. If $4M+e\le r^*\le\frac{t+R}{2}$, then $r\thicksim r^*$ and hence
\begin{equation}\label{602}
\begin{aligned}
b_q(t, r)&\lesssim \int_{0}^{1}e^{-\lambda (t+R)}e^{\lambda r^*} \lambda^{q-1}d\lambda\\
&\lesssim\int_{0}^{1}e^{-\frac{\lambda (t+R)}{2}}\lambda^{q-1}d\lambda\\
&\lesssim (t+R)^{-q}\int_0^{\infty}e^{-\theta}\theta^{q-1}d\theta\\
&\lesssim  (t+R)^{-q},
\end{aligned}
\end{equation}
which also holds for $b_{q+1}$
\begin{equation}\label{603}
\begin{aligned}
b_{q+1}(t, r)&\lesssim (t+R)^{-(q+1)}.
\end{aligned}
\end{equation}
While if $\frac{t+R}{2}\le r^*\le t+R$, we have
\begin{equation}\label{604}
\begin{aligned}
b_q&\lesssim \int_0^1\frac{\lambda^{q-1}}{(1+\lambda)r}d\lambda\\
&\lesssim \int_0^1\frac{\lambda^{q-1}}{1+\lambda(t+R)}d\lambda\\
&\lesssim (t+R)^{-q}\int_0^{\infty}(1+\theta)^{-1}\theta^{q-1}d\theta\\
&\lesssim (t+R)^{-q}
\end{aligned}
\end{equation}
and
\begin{equation}\label{605}
\begin{aligned}
b_{q+1}&\lesssim \int_0^1e^{-\lambda(t+R+1-r^*)}\frac{1}{\lambda r^*}\lambda^{q}d\lambda\\
&\lesssim (t+R)^{-1}\int_0^1e^{-\lambda(t+R+1-r^*)}\lambda^{q-1}d\lambda\\
&\lesssim (t+R)^{-1}(t+R+1-r^*)^{-q}\int_0^{\infty}e^{-\theta}\theta^{q-1}d\theta\\
&\lesssim (t+R)^{-1}(t+R+1-r^*)^{-q},\\
\end{aligned}
\end{equation}
and we finish the proof Lemma \ref{bq}.

Set
\[
\eta^*(t)=\eta_{\chi_{[\frac18, 1]}}(t),~~\eta_T^*(t)=\eta^*\left(\frac tT\right)\\
\]
and for $L\in[16R, T]\subset[16R, T(\e)]$
\begin{equation}\label{606}
\begin{aligned}
Y\left[b_q|u|^p\right](L)=\int_1^L\left(\int_0^T\int_{2M}^{\infty}
\int_{\mathbb{S}^2}|u|^pb_q(\eta_\sigma^*)^{2p'}\alpha_L^{2p'}(r^*)r^2
d\omega drdt\right)\sigma^{-1}d\sigma.
\end{aligned}
\end{equation}
For simplicity, we denote $Y(L)$ for $Y\left[b_q|u|^p\right](L)$ and have
\begin{equation}\label{0606}
\begin{aligned}
LY'(L)&=\int_0^T\int_{2M}^{\infty}
\int_{\mathbb{S}^2}|u|^pb_q(\eta_L^*)^{2p'}\alpha_L^{2p'}(r^*)r^2d\omega drdt\\
&\gtrsim\int_{\frac L{16}}^{\frac{2L}{3}}\int_{2M}^{\infty}
\int_{\mathbb{S}^2}|u|^pb_q\eta_L^{2p'}\alpha_L^{2p'}(r^*)r^2d\omega drdt\\
&\gtrsim L^{-q}\e^pL^{3-p}\\
&=\e^p,\\
\end{aligned}
\end{equation}
where the lower bounds \eqref{65}, $\eqref{66}_2$ and $p=1+\sqrt2, q=2-\sqrt2$ have been used. Also,
\begin{equation}\label{607}
\begin{aligned}
Y(L)&=\int_0^T\int_{2M}^{\infty}
\int_{\mathbb{S}^2}|u|^pb_q\alpha_L^{2p'}(r^*)r^2\left(\int_1^L
(\eta_\sigma^*)^{2p'}\sigma^{-1}d\sigma\right)d\omega drdt\\
&=\int_0^T\int_{2M}^{\infty}
\int_{\mathbb{S}^2}|u|^pb_q\alpha_L^{2p'}(r^*)r^2\left(\int_{\frac tL}^t
(\eta^*)^{2p'}(s)s^{-1}ds\right)d\omega drdt\\
&\lesssim\int_0^T\int_{2M}^{\infty}
\int_{\mathbb{S}^2}|u|^pb_q\alpha_L^{2p'}(r^*)r^2\left(\int_{\max(\frac tL, \frac18)}^{\frac23}
\eta^{2p'}(s)s^{-1}ds\right)d\omega drdt\\
&\lesssim\int_0^T\int_{2M}^{\infty}
\int_{\mathbb{S}^2}|u|^pb_q\alpha_L^{2p'}(r^*)\eta^{2p'}\left(\frac tL\right)r^2
\left(\int_{\frac18}^{\frac23}
s^{-1}ds\right)d\omega drdt\\
&\lesssim\int_0^L\int_{2M}^{\infty}
\int_{\mathbb{S}^2}|u|^pb_q\alpha_L^{2p'}(r^*)\eta_L^{2p'}(t)r^2d\omega drdt.\\
\end{aligned}
\end{equation}
Using \eqref{607}, we may bound $Y'(L)$ from below by $Y^p$. Multiplying the equation \eqref{Smain} with $b_q\alpha_L^{2p'}(r^*)\eta_L^{2p'}(t)r^2$ and then integrating with respect to $r, \omega, t$ we come to
\begin{equation}\label{608}
\begin{aligned}
&\int_0^L\int_{2M}^{\infty}
\int_{\mathbb{S}^2}|u|^pb_q\alpha_L^{2p'}(r^*)\eta_L^{2p'}(t)r^2d\omega drdt\\
\lesssim&\int_0^L\int_{2M}^{\infty}
\int_{\mathbb{S}^2}\frac{r^2}{1-\frac{2M}{r}}u\partial_t^2\left(b_q\alpha_L^
{2p'}\eta_L^{2p'}\right)\\
&-u\partial_r\left(r(r-2M)\partial_r(b_q\alpha_L^{2p'}
\eta_L^{2p'})\right)d\omega drdt\\
\lesssim&\int_0^L\int_{2M}^{\infty}
\int_{\mathbb{S}^2}\left(\frac{r^2}{1-\frac{2M}{r}}2|u\partial_tb_q\alpha_L^
{2p'}\partial_t\eta_L^{2p'}|+\frac{r^2}{1-\frac{2M}{r}}|ub_q\alpha_L^
{2p'}\partial_t^2\eta_L^{2p'}|\right)d\omega drdt\\
&+\int_0^L\int_{2M}^{\infty}
\int_{\mathbb{S}^2}|u|\Big(|2(r-M)b_q\partial_r\alpha_L^{2p'}\eta_L^{2p'}|\\
&+|2r(r-2M)\partial_rb_q\partial_r\alpha_L^{2p'}\eta_L^{2p'}|
+|r(r-2M)b_q\partial_{rr}\alpha_L^{2p'}\eta_L^{2p'}|\Big)d\omega drdt\\
\lesssim& I_{12}+I_{13}+I_{14}+I_{15}+I_{16}.\\
\end{aligned}
\end{equation}
All the five terms $I_{12}-I_{16}$ which remain to be estimated include at least one derivative over the cut-off functions $\eta_L$ or $\alpha_L$, which will restrict the time variable over $[\frac L3, \frac{2L}{3}]$ and
$[\frac L{16}, \frac{2L}{3}]$ respectively. The latter restriction with left endpoint $\frac L{16}$ is due to
\[
\frac L{16}\le r^*\le t+R.
\]
We estimate $I_{12}, I_{13}$ first, which read by combining the asymptotic behavior $\eqref{66}_2, \eqref{66}_3$ and the fact $1-\frac{2M}{r}\thicksim C$
\begin{equation}\label{609}
\begin{aligned}
I_{12}\triangleq&\int_0^L\int_{2M}^{\infty}
\int_{\mathbb{S}^2}\frac{r^2}{1-\frac{2M}{r}}2|u\partial_tb_q\alpha_L^
{2p'}\partial_t\eta_L^{2p'}|d\omega drdt\\
\lesssim&L^{-1}\left(\int_0^L\int_{2M}^{\infty}
\int_{\mathbb{S}^2}|u|^pb_q(\eta_L^*)^{2p'}\alpha_L^{2p'}r^2d\omega drdt\right)^{\frac1p}\\
&\times\left(\int_{\frac L3}^{\frac{2L}{3}}\int_{4M+e}^{t+R}b_{q}^{-\frac{1}{p-1}}b_{q+1}^{\frac{p}{p+1}}
r^2\left(1-\frac{2M}{r}\right)^{1-p'}dr^*dt
\right)^{\frac1{p'}}\\
&\lesssim (\ln L)^{\frac{p-1}{p}}\left(\int_0^L\int_{2M}^{\infty}
\int_{\mathbb{S}^2}|u|^pb_q(\eta_L^*)^{2p'}\alpha_L^{2p'}r^2d\omega drdt\right)^{\frac1p},\\
\end{aligned}
\end{equation}
\begin{equation}\label{610}
\begin{aligned}
I_{13}\triangleq&\int_0^L\int_{2M}^{\infty}
\int_{\mathbb{S}^2}\frac{r^2}{1-\frac{2M}{r}}|ub_q\alpha_L^
{2p'}\partial_t^2\eta_L^{2p'}|d\omega drdt\\
&\lesssim \left(\int_0^L\int_{2M}^{\infty}
\int_{\mathbb{S}^2}|u|^pb_q(\eta_L^*)^{2p'}\alpha_L^{2p'}r^2d\omega drdt\right)^{\frac1p}.\\
\end{aligned}
\end{equation}
Next we move to the term $I_{15}$ and have by $\eqref{66}_1, \eqref{66}_3$
\begin{equation}\label{611}
\begin{aligned}
I_{15}\triangleq&\int_0^L\int_{2M}^{\infty}
\int_{\mathbb{S}^2}|u||2r(r-2M)\partial_rb_q\partial_r\alpha_L^{2p'}\eta_L^{2p'}|
d\omega drdt\\
\lesssim&L^{-1}\left(\int_0^L\int_{2M}^{\infty}
\int_{\mathbb{S}^2}|u|^pb_q(\eta_L^*)^{2p'}\alpha_L^{2p'}r^2d\omega drdt\right)^{\frac1p}\\
&\times\left(\int_{\frac L8}^{\frac{2L}{3}}\int_{2M}^{\infty}
(r(r-2M))^{p'}b_q^{-\frac{1}{p-1}}b_{q+1}^{\frac{p}{p-1}}
r^{-\frac{2}{p-1}}\left(1-\frac{2M}{r}\right)dr^*dt\right)^{\frac1{p'}}\\
\lesssim&(\ln L)^{\frac{p-1}{p}}\left(\int_0^L\int_{2M}^{\infty}
\int_{\mathbb{S}^2}|u|^pb_q(\eta_L^*)^{2p'}\alpha_L^{2p'}r^2d\omega drdt\right)^{\frac1p},\\
\end{aligned}
\end{equation}
while $I_{14}, I_{16}$ can be estimated in the same way to get
\begin{equation}\label{612}
\begin{aligned}
I_{14}, I_{16}
&\lesssim \left(\int_0^L\int_{2M}^{\infty}
\int_{\mathbb{S}^2}|u|^pb_q(\eta_L^*)^{2p'}\alpha_L^{2p'}r^2d\omega drdt\right)^{\frac1p}.\\
\end{aligned}
\end{equation}
Finally we conclude from the definition \eqref{606}, \eqref{0606}, $\eqref{607}-\eqref{612}$ that
\begin{equation}\label{613}
\left\{
\begin{array}{ll}
&LY'(L)\gtrsim\e^p,\\
&L(\ln L)^{p-1}Y'(L)\gtrsim Y^p(L).\\
\end{array}
\right.
\end{equation}
We then can apply the following lemma with $p_1=p_2=p=1+\sqrt2$
and $\delta=\e^p$ to system \eqref{613} to get the upper bound of lifespan estimate \eqref{lf2} in Theorem \ref{thm2}, due to the fact $L$ is arbitrary in $[16R, T)$.
\begin{lem}
{\bf (Lemma 3.10 in \cite{ISWa})}. Let $2<t_0<T$. $0\le \phi\in C^1([t_0, T))$. Assume that
\begin{equation}
\left\{
\begin{aligned}
& \delta\le K_1t\phi'(t), \quad t\in (t_0, T), \\
& \phi(t)^{p_1}\le K_2t(\log t)^{p_2-1}\phi'(t), \quad t\in (t_0, T)\\
\end{aligned}
\right.
\end{equation}
with $\delta, K_1, K_2>0$ and $p_1, p_2>1$. If $p_2<p_1+1$, then there exists positive constants $\delta_0$ and $K_3$(independent of $\delta$) such that
\begin{equation}
\begin{aligned}
T\le \exp\left(K_3\delta^{-\frac{p_1-1}{p_1-p_2+1}}\right)
\end{aligned}
\end{equation}
when $0<\delta<\delta_0$.
\end{lem}
\begin{rem}
Once the differential system \eqref{613} is established, one can also use a direct method as that in the end of \cite{LLTW} to get the desired upper bound of lifespan estimate.
\end{rem}

\section{Acknowledgement}

N. A. Lai would like to express his sincere thank to Dr. Masahiro Ikeda for the helpful discussion for the critical case, and to Dr. Kyouhei Wakasa for the long time communication on this problem.

The first author was supported by NSFC (No. 12271487), the second author was supported by NSFC (No. 12171097).

\bibliography{mybibfile}

\begin{thebibliography}{99}




\bibitem{BS}
P. Blue and J. Sterbenz, Uniform decay of local energy and the semi-linear wave equation on Schwarzschild space, Comm.
Math. Phys., 268 (2006) 481-504.

\bibitem{Catania}
D. Catania and V. Georgiev, Blow-up for the semilinear wave equation in the Schwarzschild metric, Differential Integral Equations, 19 (2006) 799-830.
\bibitem{DR}
M. Dafermos and I. Rodnianski, Small-amplitude nonlinear waves on a black hole
background, J. Math. Pures Appl., 84 (2005) 1147-1172.

\bibitem{DFW}
W. Dai, Y. D. Fang and B. C. Wang,
Lifespan of solutions to the Strauss type wave system on asymptotically flat space-times, Discrete Contin. Dyn. Syst.
Ser. S., 40(8) (2020) 4985-4999.



\bibitem{Georgive}
V. Georgiev, H. Lindblad, C. D. Sogge, Weighted Strichartz estimates and global existence for semilinear wave equations,
Amer. J. Math., 119(6) (1997) 1291-1319.

\bibitem{Glassey1}
R. T. Glassey, Finite-time blow-up for solutions of nonlinear wave equations, Math. Z., 177 (1981) 323-340.
\bibitem{Glassey2}
R. T. Glassey, Existence in the large for $\Box u=F(u)$ in two space dimensions, Math. Z., 178 (1981) 233-261.

\bibitem{ISWa}M.Ikeda, M.Sobajima and K. Wakasa,
Blow-up phenomena of semilinear wave equations and their weakly coupled systems, J. Differential Equations, 267(9) (2019) 5165-5201.


\bibitem{Takamura6}
T. Imai, M. Kato, H. Takamura and K. Wakasa,
%The sharp lower bound of the lifespan of solutions to semilinear wave equations with low powers in two space dimensions,  Adv. Stud. Pure Math. 81, 31--53 (2019).
The sharp lower bound of the lifespan of solutions to semilinear
  wave equations with low powers in two space dimensions.
In {\em {Asymptotic analysis for nonlinear dispersive and wave
  equations. Proceedings of the international conference on asymptotic analysis
  for nonlinear dispersive and wave equations, Osaka University, Osaka, Japan,
  September 6--9, 2014}}, pages 31--53. Tokyo: Mathematical Society of Japan,
  2019.

\bibitem{John}
F. John, Blow-up of solutions of nonlinear wave equations in three space dimensions,
Manuscripta Math., 28 (1-3)(1979) 235-268.

\bibitem{LT1}
N. A. Lai and Z. H. Tu, Strauss exponent for semilinear wave equations with scattering space dependent damping, J. Math. Anal. Appl., 489 (2020) 124189.

 \bibitem{LLWW}%[53]
N. A. Lai, M. Y. Liu, K. Wakasa and C. B. Wang. Lifespan estimates for
2-dimensional semilinear wave equations in asymptotically Euclidean exterior domains, J. Funct. Anal., 281 (12)(2021) 109253.

\bibitem{LLTW}
N. A. Lai, M. Y. Liu, Z. H. Tu and C. B. Wang. Lifespan estimates for semilinear wave equations with space dependent damping and potential, arXiv:2102.10257.


\bibitem{LZ}
N. A. Lai and Y. Zhou, An elementary proof of Strauss conjecture, J. Funct. Anal., 267 (5)(2014) 1364-1381.


\bibitem{LaiZ}
N. A. Lai and Y. Zhou, Lifespan estimate for the semilinear wave equation with a derivative nonlinear term in
Schwarzschild spacetime (in Chinese), Sci. Sin. Math., 51 (2021) 957-970. doi: 10.1360/SSM-2020-0233.
%\bibitem{LCM}
%N. A. Lai, L. Chen and Z. Y. Ma,
%Finite time blow up for a wave equation with nonlinear memory in Schwarzschild
%spacetime(in Chinese), Sci. Sin. Math., 45 (2015), 117-128.

\bibitem{LZbook}
T. T. Li and Y. Zhou, Nonlinear Wave Equations(in Chinese), Series in Contemporary Mathematics, Shanghai Scientific \& Technical Publishers, 2016.

\bibitem{LLM}
Y. H. Lin, N. A. Lai and S. Ming, Lifespan estimate for semilinear wave equation in Schwarzschild spacetime, Appl Math Lett., 99 (2020) 105997, 4 pp. doi: 10.1016/j.aml.2019.105997.


\bibitem{Lindblad3}
H. Lindblad, Blow-up for solutions of $\Box u=|u|^p$ with small initial data, Comm. Partial Differential Equations, 15 (6) (1990)
757-821.

\bibitem{Lindblad2}H. Lindblad, C. D. Sogge, Long-time existence for small amplitude semilinear wave equations,
Amer. J. Math., 118 (5) (1996) 1047-1135.


\bibitem{LMS}
H. Lindblad, J. Metcalfe, C. D. Sogge, M. Tohaneanu and C. B. Wang, The Strauss conjecture on Kerr black hole backgrounds, Math. Ann., 359(3-4) (2014) 637-661.


\bibitem{Luk}%19-26
J. Luk, The null condition and global existence for nonlinear wave equations on slowly rotating Kerr spacetimes, Journal Eur Math Soc., 15 (2013) 1629-1700.



\bibitem{MMT}
J. Marzuola, J. Metcalfe, D. Tataru and M. Tohaneanu, Strichartz estimates on Schwarzschild black hole backgrounds, Comm. Math.
Phys., 293 (2010) 37-83.




\bibitem{MW}
J. Metcalfe and C. B. Wang, The Strauss conjecture on asymptotically flat space-times, SIAM J Math. Anal., 49 (2017) 4579-4594.


\bibitem{Nico1} J. P. Nicolas, Nonlinear Klein-Gordon equation on Schwarzschild-like metrics, J. Math.
Pures Appl., 74(9) (1995) 35-58.


\bibitem{Nico2} J. P. Nicolas, A nonlinear Klein-Gordon equation on Kerr metrics, J. Math. Pures Appl., 81(9) (2002) 885-914.



\bibitem{Schaeffer}
J. Schaeffer, The equation $\Box u=|u|^p$ for the critical value of $p$, Proc. Roy. Soc. Edinburgh
Sect. A, 101 (1-2)(1985) 31-44.
\bibitem{Sideris}T. C. Sideris, Nonexistence of global solutions to semilinear wave equations in high dimensions, J. Differential Equations, 52 (1984) 378-406.


\bibitem{Takamura2}
  H. Takamura, Improved Kato's lemma on ordinary differential inequality and its
application to semilinear wave equations, Nonlinear Analysis, 125 (2015) 227-240.

\bibitem{Takamura}
  H. Takamura and K. Wakasa, The sharp upper bound of the lifespan of solutions to critical semilinear wave equations in high
dimensions, J. Differential Equations, 251 (2011) 1157-1171.



\bibitem{Yordanov}
B. Yordanov and Q. S. Zhang, Finite time blow up for critical wave equations in high dimensions, J. Funct. Anal., 231 (2006) 361-374.

\bibitem{Wan1}
C. B. Wang, Long-time existence for semilinear wave equations on asymptotically flat space-times, Comm. Partial
Differential Equations, 42 (2017) 1150-1174.

\bibitem{Wan}
C. B. Wang,
Recent progress on the Strauss conjecture and related problems(in Chinese), Sci. Sin. Math., 48 (2018)
111-130.

\bibitem{Zhou3}
Y. Zhou, Lifespan of classical solutions to $\Box u=|u|^p$ in two space dimensions, Chin. Ann. Math. Ser.B, 14 (1993) 225-236.

\bibitem{Zhou5}
Y. Zhou, Blow up of classical solutions to $\Box u=|u|^{1+\alpha}$ in three space dimensions, J. Partial Differential Equations, 5 (1992)
21-32.

\bibitem{Zhou2}
Y. Zhou, Blow up of solutions to semilinear wave equations with critical exponent in high
dimensions, Chinese Ann. Math. Ser. B, 28(2)(2007) 205-212.

\end{thebibliography}

\bibliographystyle{plain}
%\bibliography{bibli}
\end{CJK*}

\end{document}